\documentclass{amsart}
\usepackage{hyperref}

\theoremstyle{definition}

\theoremstyle{definition} 
\newtheorem*{remark*}{Remark}

\newcommand{\dd}{\operatorname{d}}
\renewcommand{\dd}{\,\mathrm{d}}

\newcommand{\ii}{\operatorname{I}}

\newcommand{\vp}{\varepsilon}

\begin{document}

\title{On the H\"older and Cauchy--Schwarz inequalities}
\markright{H\"older and Cauchy--Schwarz}
\author{Iosif Pinelis}
\address{Department of Mathematical Sciences, 
Michigan Technological University, 
Houghton, Michigan 49931\\
ipinelis@mtu.edu}

\begin{abstract}
A generalization of the H\"older inequality is considered. Its relations with a previously obtained improvement of the Cauchy--Schwarz inequality are discussed. 
\end{abstract}

\maketitle

Let $f$ and $g$ be any nonnegative measurable functions on a measure space $(S,\Sigma,\mu)$. Let 
$p$ be be any real number in the interval $(1,\infty)$, and let then $q:=\frac p{p-1}$, so that 
$\frac1p+\frac1q=1$. H\"older's inequality states that 
\begin{equation}\label{eq:holder}
	\mu(fg)\le\mu(f^p)^{1/p}\mu(g^q)^{1/q}, 
\end{equation}
where $\mu(f):=\int_S f\dd\mu\in[0,\infty]$. 

Consider any Borel-measurable transformation 
\begin{equation*}
	[0,\infty)^2\ni(x,y)\longmapsto T(x,y)=\big(T_1(x,y),T_2(x,y))\in[0,\infty)^2 
\end{equation*}
preserving the product: 
\begin{equation*}
	T_1(x,y)\,T_2(x,y)=xy
\end{equation*}
for all $(x,y)\in[0,\infty)^2$. Examples of such product-preserving transformations are given by the formulas $T(x,y)=(kx,y/k)$ for positive real $k$, $T(x,y)=(y,x)$, $T(x,y)=(x\vee y,x\wedge y)$, and any compositions thereof. 

For each $j\in\{1,2\}$, let $T_j(f,g)$ be the nonnegative Borel-measurable function on $S$ defined by the ``pointwise'' formula $T_j(f,g)(s):=T_j\big(f(s),g(s)\big)$ for all $s\in S$. 
Then trivially $T_1(f,g)T_2(f,g)=fg$. So, substituting $T_1(f,g)$ and $T_2(f,g)$ for $f$ and $g$ in  \eqref{eq:holder}, one immediately has the following generalization of H\"older's inequality: 
\begin{equation}\label{eq:holderT}
	\mu(fg)\le 
	\mu\big(T_1(f,g)^p\big)^{1/p}\,\mu\big(T_2(f,g)^q\big)^{1/q}. 
\end{equation}
In particular, choosing the product-preserving transformation defined by the formula 
$T(x,y)=(x\vee y,x\wedge y)$\ for all\ $(x,y)\in[0,\infty)^2$, one has 
\begin{equation}\label{eq:maxmin_p}
	\mu(fg)\le B_p(f,g):=
	\mu\big((f\vee g)^p\big)^{1/p}\,\mu\big((f\wedge g)^q\big)^{1/q}. 
\end{equation}

In the Cauchy--Schwarz case, when $p=q=2$, inequality \eqref{eq:maxmin_p} can be rewritten as  
\begin{equation}\label{eq:c-s}
	\mu(fg)^2\le\mu(a\vee b)\,\mu(a\wedge b).  
\end{equation}
Here and in what follows, set 
\begin{equation*}
	a:=f^2\quad\text{and}\quad b:=g^2. 
\end{equation*}
Under the additional assumption that the measure space is the Lebesgue measure space over an interval, inequality \eqref{eq:c-s} was given in \cite{xiang13}. The proof in \cite{xiang13} is rather complicated and uses a discrete approximation; it also appears to critically depend on the Cauchy--Schwarz assumption $p=2$.  

On the other hand, it was noted in \cite{xiang13} that the upper bound in \eqref{eq:c-s} improves the Cauchy--Schwarz bound $\mu(a)\,\mu(b)$ on $\mu(fg)^2$. Indeed, one has the identity 
\begin{equation}\label{eq:improv}
	\mu(a)\,\mu(b)=\mu(a\vee b)\,\mu(a\wedge b)+\mu\big((a-b)_+\big)\,\mu\big((b-a)_+\big), 
\end{equation}
where $u_+:=u\vee0$. 
To quickly verify this identity, note that 
\begin{align*}
	\mu(a)\mu(b) =& \mu\big(a \vee b - (b - a)_+ \big)\mu\big(a \wedge b + (b - a)_+ \big) \\ 
=& \big[\mu(a \vee b) - \mu\big((b - a)_+ \big)\big]\big[\mu(a \wedge b) + \mu\big((b - a)_+ \big)\big] \\ 
=& \mu(a \vee b)\mu(a \wedge b) + \mu(a \vee b)\mu\big((b - a)_+ \big) \\ 
&- \mu\big((b - a)_+ \big)\mu(a \wedge b) - \mu\big((b - a)_+ \big)\mu\big((b - a)_+ \big) \\ 
=& \mu(a \vee b)\mu(a \wedge b) + \mu\big(a \vee b - a \wedge b - (b - a)_+ )\mu((b - a)_+ \big) \\ 
=& \mu(a \vee b)\mu(a \wedge b) + \mu\big((a - b)_+ \big)\mu\big((b - a)_+ \big). 
\end{align*} 

A more general result on improvements of the Cauchy--Schwarz inequality is due to Daykin--Eliezer--Carlitz \cite{dec} in the case when $\mu$ is the counting measure on a finite set, which was then extended in \cite{sitnik} to the case when $\mu$ is the Lebesgue measure on an interval. Those results still appear to depend on the Cauchy--Schwarz assumption $p=2$. I am pleased to thank S.\ M.\ Sitnik for the reference to \cite{sitnik}. 

One may then wonder whether the ``max-min'' upper bound $B_p(f,g)$ in \eqref{eq:maxmin_p} improves the H\"older bound for all $p\in(1,\infty)$.  
However, it turns out that the Cauchy--Schwarz case $p=2$ is the only exception here. 
Even the smaller, symmetrized upper bound $B_p(f,g)\wedge B_q(f,g)$ does not in general improve the H\"older bound for any $p\in(1,\infty)\setminus\{2\}$. 
The key observation here is that, according to \eqref{eq:improv}, the improvement \break $\mu\big((a-b)_+\big)\,\mu\big((b-a)_+\big)$ of the bound in \eqref{eq:c-s} over the Cauchy-Schwarz bound is no greater than $\mu(1)^2\vp^2$ if $|a-b|\le\vp$, for any positive real $\vp$. 
So, to show that the Cauchy--Schwarz case $p=2$ is indeed exceptional, one may try to choose the functions $f$ and $g$ to be close to each other, at least when $p$ is close to $2$. 

It is not hard to make this idea work. 
Indeed, let e.g.\ the measure space $(S,\Sigma,\mu)$ be the Lebesgue measure space over the interval $[0,1)$. 
Fix any $p\in(1,\infty)\setminus\{2\}$ and any $m\in(0,1)$. 
Fix then any $w\in(0,\infty)$ such that $w\in(0,1)$ if $p\in(1,2)$ and $w\in(1,\infty)$ if $p\in(2,\infty)$ Thus, in any  case, $w^p-w^q>0$. 

Let then 
\begin{equation*}
\text{$g=g_{m,w}:=w\ii_{[0,m)}+\ii_{[m,1)}$ and $f=f_{m,w,t}:=(1-t)w\ii_{[0,m)}+(1+t)\ii_{[m,1)}$}	
\end{equation*}
for $t\in[0,1)$, so that $f$ and $g$ are positive measurable functions on $[0,1)$; here $\ii_A$ denotes the indicator function of a set $A$. 
Then 

\newpage
\begin{align*}
&	B_p(f,g)-
	\mu(f^p)^{1/p}\mu(g^q)^{1/q} 
	\\ 
&	
=\mu\big((f\vee g)^p\big)^{1/p}\,\mu\big((f\wedge g)^q\big)^{1/q}-
	\mu(f^p)^{1/p}\mu(g^q)^{1/q} \\ 
&	 
\begin{aligned}
	=d_1(t)
	:=&\big((1 - m) (1 + t)^p + m w^p\big)^{1/p} \big(1 - m + m (1 - t)^q w^q\big)^{1/q}  \\ 
	&-  \big((1 - m) (1 + t)^p + m (1 - t)^p w^p\big)^{1/p} (1 - m + m w^q)^{1/q}, 
\end{aligned}
\\
&	B_q(f,g)-
	\mu(f^p)^{1/p}\mu(g^q)^{1/q} \\ 
&	=\mu\big((f\vee g)^q\big)^{1/q}\,\mu\big((f\wedge g)^p\big)^{1/p}-
	\mu(f^p)^{1/p}\mu(g^q)^{1/q} \\ 	
&	
\begin{aligned}
=d_2(t)
:=& 
	 \big((1 - m) (1 + t)^q + m w^q\big)^{1/q} \big(1 - m + m (1 - t)^p w^p\big)^{1/p} \\ 
	&-  \big((1 - m) (1 + t)^p + m (1 - t)^p w^p\big)^{1/p} (1 - m + m w^q)^{1/q},  
\end{aligned}		
\end{align*}
so that $d_1(t)\wedge d_2(t)$ is the difference between the symmetrized ``max-min'' upper bound $B_p(f,g)\wedge B_q(f,g)$ and the H\"older bound $\mu(f^p)^{1/p}\mu(g^q)^{1/q}$, with $f=f_{m,w,t}$ and $g=g_{m,w}$. 
For each $j\in\{1,2\}$ one has $d_j(0)=0$ and 
$$d'_j(0)=(1 - m) m (w^p - w^q) (1 - m + m w^p)^{-1/q} (1 - m + m w^q)^{-1/p}>0,$$ 
so that 
$d_1(t)\wedge d_2(t)>0$ for all small enough positive $t$. 
Thus, for each $p\in(1,\infty)\setminus\{2\}$ we have constructed positive measurable functions $f$ and $g$ such that the symmetrized ``max-min'' upper bound $B_p(f,g)\wedge B_q(f,g)$ is strictly greater than the H\"older bound $\mu(f^p)^{1/p}\mu(g^q)^{1/q}$.

\bibliographystyle{abbrv}
\bibliography{C:/Users/ipinelis/Dropbox/mtu/bib_files/citations12.13.12}

\end{document}